\documentclass[11pt,twoside]{article}
\usepackage{latexsym}
\usepackage{amssymb,amsbsy,amsmath,amsfonts,amssymb,amscd}
\usepackage{graphicx}
\usepackage{hyperref}
\usepackage{setspace}
\usepackage{dsfont}

\usepackage{pdfsync}
\usepackage{color}
\newcommand{\commentout}[1]{}
\newcommand{\R}{\mathbb{R}}

\newcommand {\al} {\alpha}
\newcommand {\be} {\beta}
\newcommand {\e}  {\varepsilon}

\newcommand {\Chi} {{\bf \raise 2pt \hbox{$\chi$}} }
\newcommand {\f}   {\frac}
\newcommand {\p}   {\partial}
\newcommand {\ddt}   { \frac{d}{dt}}

\newcommand {\proof} {\noindent {\bf Proof}. }
\newcommand{\re}{\eqref}
\newcommand{\ol}{\overline}
\newcommand{\ul}{\underline}
\newcommand{\beq}{\begin{equation}}
\newcommand{\eeq}{\end{equation}}
\newcommand{\bea} {\begin{array}{rl}}
\newcommand{\eea} {\end{array}}
\newcommand{\bepa}{\left\{ \begin{array}{l}}
\newcommand{\eepa} {\end{array}\right.}
\newtheorem{theorem}{Theorem}[section]

\newtheorem{remark}[theorem]{Remark}

\newcommand{\qed}{{ \hfill
                       {\unskip\kern 6pt\penalty 500 \raise -2pt\hbox{\vrule\vbox to 6pt{\hrule width 6pt
                       \vfill\hrule}\vrule} \par}   }}

\title{\Large \bf Long-term analysis of phenotypically structured models}

\author{
Alexander Lorz\thanks{Sorbonne Universit\'es, UPMC Univ. Paris 06, UMR 7598, Laboratoire Jacques-Louis Lions, F-75005, Paris, France}
\thanks{CNRS, UMR 7598, Laboratoire Jacques-Louis Lions, F-75005, Paris, France
} 
\thanks{INRIA Paris-Rocquencourt, EPC Mamba} 
\footnotemark[4]
\and Beno\^ \i t Perthame\footnotemark[1] \footnotemark[2] \footnotemark[3] \thanks{Emails: lorz@ann.jussieu.fr, benoit.perthame@upmc.fr}
}

\date{\today}

\begin{document}
\maketitle
\pagestyle{plain}
\pagenumbering{arabic}

\begin{abstract}
Phenotypically structured equations arise in population biology to describe the interaction of species with their environment that brings the nutrients. This interaction usually leads to selection of the fittest individuals.  Models used in this area are highly nonlinear, and the question of long term behaviour is usually not solved. However, there is a particular class of models for which convergence to an Evolutionary Stable Distribution is proved, namely when the quasi-static assumption is made. This means that the environment, and thus the nutrient supply, reacts immediately to the population dynamics. One possible proof is based on a Total Variation bound for the appropriate quantity.

We extend this proof to several cases where the nutrient is regenerated with delay. A simple example is the chemostat with a rendering factor, then our result does not use any smallness assumption. For a more general setting, we can treat  the case with a fast  reaction of nutrient supply to the population dynamics.

\end{abstract}

\bigskip

 \\[4mm]


\section{Introduction}
\label{sec:intro}

In population biology, long-term behaviour for phenotypically structured models is a difficult question related to interaction with environmental conditions, selection of fittest trait and lack of dissipation principles.  The {\em competitive exclusion principle} is a famous  general result, and states that, with a single type of `niche' or substrate, a single trait is selected. 

A typical example where this can be proved rigorously, is the chemostat model 
\begin{align} 
& \f{\partial}{\partial t} n(t,x) = n \big[ -R_0 +  a(x)\eta \big(x,S(t) \big) \big],  \qquad x \in \R^d, \; t \geq 0,  \label{eq.ch_n}
\\[5pt]
& \frac{d}{dt} S(t) + R_0 S(t) = R_0 S_0 -  \int n \eta \big(x,S(t) \big) \,dx. \label{eq.ch_S} 
\end{align}
The first equations describes the population density $n(t,x)$ of individuals which at time $t$ have the trait $x$. The substrate, whose concentration is denoted by $S$, is used with a trait-dependent uptake coefficient  $\eta(x,S)$ and a rendering factor $a(x)$. The renewal of the reactor, with fresh nutrient $S_0$, occurs with the rate $R_0$. 

The simplest situation is when there is a unique Evolutionary Stable Distribution (ESD in short, a term coined in  \cite{PJ.GR:09}) which concentrates in a single Dirac mass.  That means there is a unique trait $\bar x$, associated with a nutrient concentration $\bar S >0 $, characterized by
\beq
\max_x  a(x)\eta(x,\bar S)    = R_0 =  a(\bar x)\eta(\bar x,\bar S).
\label{ch:ESS}
\eeq
The first equality allows to compute a unique $\bar S$, assuming $\eta$ is increasing with $S$. And the second equality gives $\bar x$. 

Then, it is known  when $a\equiv 1$, see \cite{BP:07}, that the competitive exclusion principle can be expressed as 
\beq
n(t,x)  \underset{t \to \infty}{\longrightarrow} \bar \rho \delta(x- \bar x),
\label{concentration}
\eeq
 and we extend this result here.

However, we do not know  general  assumptions on $\eta^n$, $\eta^S$ which would lead to a similar result for the more general chemostat model
\begin{align*} 
& \f{\partial}{\partial t} n = n \big(-R_0 +  \eta^n(x,S)\big),
\\
& \frac{d}{dt} S + R_0 S = R_0 S_0 - S \int n \eta^S(x,S) \,dx.  
\end{align*}
A general method is to use a Lyapunov functional (entropy) but this requires a particular structure on the system, \cite{PJ.GR:09,NC.PJ.GR:11,SM.BP.JW:12}. 

The laws for nutrient delivery and consumption may differ for other models \cite{RS.JE:02,chemostatbook},  but similar questions still arise. A `generic' mathematical model, which contains \eqref{eq.ch_n}--\eqref{eq.ch_S}  as a particular case,  can be written as 
\begin{align} 
& \partial_t n(t,x) = n \,R\big(x,S(t) \big), \qquad x \in \R^d, \; t \geq 0, \label{e.n} \\ 
& \beta \frac{d}{dt} S(t)  = Q\big(S(t),\rho(t) \big) ,   \label{e.I}
\\
& \rho(t) :=  \int n(t,x)\,dx.  \label{e.rho}
\end{align}
Here $R(x,S)$ denotes a generic trait-dependent birth-death rate, $S$ is still the nutrient concentration and $\rho(t)$ a measure of the pressure exerted by the total population for nutrient consumption with rate $Q$. The parameter $\beta$, which obviously could be included in $Q$ is used here for a simple mathematical purpose. It gives a time scale which, in the limit $\beta =0$, just gives $ 0 = Q\big(S(t),\rho(t) \big) $. Under suitable assumptions, this equation can be inverted in $S= q(\rho)$. In this case,  the long term selection of the ESD,   \eqref{concentration}, is known to hold \cite{GB.BP:07,AL.SM.BP:10,ackleh_F_T}.

Our aim is to prove the same convergence result  to an ESD, \eqref{concentration}, when  $\beta $ is small. Section~\ref{sec:gen} is devoted to prove the result and a precise statement is given in the Theorem \ref{th:betasmall}. In order to make the proof more intuitive, we begin with the simpler case of the chemostat system \eqref{eq.ch_n}--\eqref{eq.ch_S}; this is developed in Section \ref{sec:chemostat}.

\section{The chemostat with rendering factor}
\label{sec:chemostat}

The model  of the chemostat with a rendering factor is defined by the system \eqref{eq.ch_n}--\eqref{eq.ch_S}.  We complete it with initial data $n^0(x)$, $S^0$ that satisfy
\beq
0< S^0 \leq S_0, \qquad n^0(x)  >0   \quad \forall x \in \R^d, \quad n^0 \in L^1(\R^d).
\label{e.id} 
\eeq
We recall that the notation
$$
\rho(t)= \int_{\R^d}Ên(t,x) dx.
$$

In order to analyze the long term behaviour, we need  assumptions on the problem parameters and coefficients. Namely, we need  to ensure first non-extinction which follows from the assumptions
\beq
 \eta(x,S_0) > R_0, \quad \eta(x,0)=0, \qquad \forall x \in \R^d.
\label{e.s0} 
\eeq
Next, it is intuitive to assume that, the more nutrient available, the higher the growth rate 
\beq
0< \ul K_\eta \leq  \eta_S(x,S) \le \ol K_\eta , \qquad \forall x \in \R^d, \; \forall S\in (0, S_0).
\label{e.eta}  
\eeq
For the rendering factor, based on a biological interpretation, it is usually assumed that $ a(x) \leq 1$ but here we only use that for some constants $a_m, \; a_M >0$
\beq
0< a_m< a(x) \leq a_M. 
\label{e.a} 
\eeq

Then, we have the following generalization of the case $a\equiv 1$ which is treated in \cite{BP:07}.
\begin{theorem}
With assumptions \re{e.id}--\re{e.a}, there are constants $\rho_m, \; \rho_M$, such that 
$$
0<\rho_m < \rho(t) \le \rho_M,  \qquad 0< S(t) \leq S_0.
$$
Assuming also \eqref{ch:ESS}, as $t \to \infty$, 
$$
S(t) \to \bar S, \qquad \rho(t)  \to \bar \rho >0, \qquad  n(t,x) \rightharpoonup \bar \rho \delta(x - \bar{x}).
$$
\label{th:chemostat}
\end{theorem}

\proof {\em 1st Step. A conserved quantity.} For future use, we define 
\beq
u(t) = \int \f{n(t,x)}{a(x)}\,dx + S(t) -S_0, \qquad J :=  \ddt \int \f{n(t,x)}{a(x)}\,dx.
\label{ch_def}
\eeq
Dividing  equation \re{eq.ch_n} by $a$, integrating and adding equation  \re{eq.ch_S}, we obtain 
\[
 \frac{d}{dt}u(t) + R_0 u(t) = 0.
\]
It follows that
\beq 
u(t) =   u(0) e^{-R_0 t }, \qquad  \ddt u(t) = -R_0 u(0) e^{-R_0 t }= J + \f{dS}{dt}  .
\label{ch_exp}
\eeq

The a priori bounds follow easily. Because $n > 0$, we find $S\leq S_0$ and because $\eta(x, 0)=0$ from assumption \eqref{e.s0}, we find $S> 0$. For an upper bound on $\rho$, we use that $u(t)$ is bounded and assumption \eqref{e.a}. We find
$$
\f{\rho(t)}{a_M} \leq \int \f{n(t,x)}{a(x)}\,dx \leq \max u(t) + S_0.
$$
The lower bound $\rho_m$ can be derived in the same way, using $a_m$.
\\

\noindent {\em 2nd Step. BV Estimates of $\int \f{n(t,x)}{a(x)} dx$.}
%
Then, we can apply the argument in \cite{BP:07} which we recall now. Using the definition of $J$ in \eqref{ch_def}, we have, using \eqref{ch_exp}, 
\begin{align*}
\ddt {J} &=  \int \f{n}{a}\big( -R_0 + a(x) \eta(x,S(t)) \big)^2\,dx + \f{dS}{dt} \int n \eta_S(x, S(t))  \,dx 
\\
&\geq \f{dS}{dt} \int n \eta_S(x, S(t))  \,dx  
\\
&= \left( -R_0 u(0) e^{-R_0 t} - J\right) \int {n} \eta_S\,dx .
\end{align*}

We define the negative part of $J$ by $J_- = \max(0,-J)$. Then, we obtain 
$$
\ddt  {J}_-  +  J_- \int {n} \eta_S\,dx \leq R_0 |u(0)| \rho_M \ol K_\eta e^{-R_0 t}, 
$$
$$
\ddt  {J}_-  + \rho_m  \ul K_\eta J_-  \leq R_0 |u(0)| \rho_M \ol K_\eta e^{-R_0 t}.
$$
This proves that $J_- (t) \leq J_- (0) e^{-\nu t}$ with $\nu = \min (R_0,  \rho_m  \ul K_\eta)$. Therefore $J_- \in L^1(0, \infty)$, and because $J$ is bounded, we obtain that $J \in L^1(0, \infty)$. Therefore, $J$ has bounded variations and $\lim_{t\to \infty} \int \f{n(t,x)}{a(x)}\,dx$ exists. Because $u(t)$ converges to $0$, we conclude that $S(t)$ has a limit
$$
S(t) \underset{t \to \infty}{\longrightarrow} S_\infty.
$$ 
\\
\noindent {\em 3rd Step. The limits.}
%
At this stage we can identify $S_\infty$. This is done with the usual arguments in the field \cite{OD.PJ.SM.BP:05,LD.PJ.SM.GR:08}. From the equation \eqref{eq.ch_n}, and the bounds on $\rho$, we immediately conclude that the growth rate should vanish on the long term, that is written 
$\max_x [R_0 -a(x) \eta(x, S_\infty) ]=0$. By monotony in $S$ of $\eta$, this tells us that $S_\infty = \ol S$ and that $n(t,x)$ concentrates as a Dirac mass at the point $\bar x$ where this maximum is achieved.  This identifies completely the limits. From the limit of $S(t)$ and $u(t)$, we know that $\int \f{n(t,x)}{a(x)}dx$ converges to $S_0 - \ol S$. And from the concentration at $\bar x$,  we conclude that $ \rho(t) = \int {n(t,x)} dx$ converges to $a(\bar x) (S_0 - \ol S)$.

The Theorem \ref{th:chemostat} is proved. 
\qed

\section{The general setting}
\label{sec:gen}

In the general setting of the system \eqref{e.n}--\eqref{e.rho}, the same proof does not apply per se. This is because we do not dispose of a quantity, as $u(t)$ in the previous proof, which is easy to control and brings us back to the quasi-steady state where $S$ is a function of $\rho$. For this reason, we need a smallness condition which is well expressed in terms of $\beta$. With this condition, we can build a quantity which belongs to $BV(0,\infty)$, as $J(t)$ in the previous proof. 

\subsection{Assumptions and main result}

We complete the system \eqref{e.n}--\eqref{e.rho} with initial conditions $S^0$, $n^0$,  which are compatible with some invariant region of interest
\beq
S_m < S^0 < S_0, \qquad n^0(x)>0, \; \forall x \in \R^d, \quad n^0 \in L^1(\R^d),
\label{gen:id}
\eeq
(see the definition of $S_m$ and $S_0$ in the assumptions below, this assumption for $S^0$ is made to simplify the statements and can be seen as a generalization of those for the chemostat in Section~\ref{sec:chemostat}).

Next, for the  Lipschitz continuous functions $R$ and $Q$, we assume that there are constants $S_0>0$, $K_Q>0$... such that 
\begin{align}
& Q(0, \rho) > 0, \quad Q(S_0, \rho) \leq 0,  \; \forall \rho \geq 0, \qquad Q_S(S, \rho) \leq - K_Q , \quad   Q_\rho (S, \rho)  \le -K_Q, \label{e.Qbounds}
\\
&0<\ul K_1 \leq R_S(x,S) \leq \ol K_1  \label{e.R_I},
\\
&\sup_{0\leq S \leq S_0} \| R(\cdot,S) \|_{W^{2,\infty}(\R^d)} \leq K_2 \label{e.Rsup}.
\end{align}

Note that from assumption \eqref{e.Qbounds}, we directly obtain the bounds
\beq
n(t,x) >0, \qquad 0< S(t) \leq S_0.
\label{gen:bounds}
\eeq

With these assumptions, the smallness condition on $\beta$ can be written as 
\beq
 \min_{\begin{array}{l} 0\leq \rho \leq \rho_M, \\ S_m \leq  S \leq S_0\end{array}} \f{|Q_S|}{ |Q_\rho|} 
 \geq  \; 4 \beta \max_{\begin{array}{l} 0\leq \rho \leq \rho_M, \\ S_m \leq S \leq S_0\end{array}} \f{ \ol K_1 \rho_M}{|Q_S|}
\label{e.beta_rho}
\eeq
(see the definition of $\rho_M$, $S_m$ below, which only depends on the assumptions above).

\begin{theorem}
With assumptions \re{gen:id}--\re{e.Rsup}, there is a constant $\rho_M$ such that
$$ 
\rho(t) \le \rho_M, \quad \text{and} \quad  S_m \leq S(t), 
$$
where the value $S_m < S_0$ is defined by $Q(S_m,\rho_M)=0$ and this exists thanks to the assumption \eqref{e.Qbounds}.

Assuming also \eqref{e.beta_rho}, $\rho^2(t)$ has bounded total variation. Consequently, there are limits  $0\leq \ol \rho \leq \rho_M $, $S_m \leq \ol S \leq S_0$ 
$$
S(t)  {\longrightarrow}  \bar S , \qquad \rho(t)  {\longrightarrow}  \bar \rho, \quad \text { as }Ê\; t \to \infty ,
$$
and
$$
Q( \bar \rho, \bar S)=0, \qquad R(x,\bar S) \leq 0  \quad \forall x \in \R^d.
$$
\label{th:betasmall}
\end{theorem}

As for the chemostat, the solution can go extinct, that means $\ol \rho =0$. When  $\ol \rho >0$, from the usual methodology developed in \cite{OD.PJ.SM.BP:05, GB.BP:07, LD.PJ.SM.GR:08}, we can also conclude that
$$
0 = \max_{x} R(x, \bar  S).
$$
And,  the population density $n(t)$ concentrates on the maximum points of $R(\cdot, \bar S)$.  For instance, with the additional assumption that there is a single $ \bar{x} \in \R^d$ such that 
\beq
R(\bar x, \bar S)=0 = \max_{x} R(x, \bar  S) ,    \label{e.Rmin} 
\eeq
we  have, in the sense of measures,  
$$
n(t,x) \longrightarrow \bar \rho \delta(x-\bar{x}),
$$
that is a monomorphic population in the language of adaptive dynamics \cite{OD:04,RC.JH:05,SG.EK.GM.JM:98}.

The end of this section is devoted to prove Theorem~\ref{th:betasmall}. This requires to adapt the method introduced in \cite{BP:07,GB.BP:07, GB.BP:08} which is to prove that $\rho(t)$ has a bounded Total Variation. This method works well in the quasi-static case, that is $\beta =0$. The adaptation is not as direct as  one could think in view of Section~\ref{sec:chemostat}.

\subsection{An upper bound for $\rho$ }

This step is not as simple as usual. Integrating the equation \eqref{e.n} with respect to $x$, yields that
\[  
\ddt\rho \leq \rho ( K_2 + \ol K_1 S ), \qquad  \ddt \ln \rho  \leq K_2 + \ol K_1 S_0.
\]
Because, from our assumptions on $Q$, there are constants such that $Q(S,\rho) \le -K_3 \rho + K_4$, adding equation \re{e.I} we  obtain the inequality
\begin{align*}
\ddt(\ln \rho + \beta S) &\leq K_2 + \ol K_1 S_0 + K_4 - K_3 \rho  \leq  K_2 + \ol K_1 S_0 + K_4 - \f{ K_3}{e^{\beta S_0}} e^{ \ln \rho +  \beta S}.
\end{align*}
Therefore, for $C_2$ the root in $\ln \rho + \beta S$ of the right hand side, we have the bound
$$
 \ln \rho \leq \ln \rho + \beta S \leq \max ( \ln\rho_0 + \beta S_0, \, C_2). 
$$
This directly gives an upper bound $\rho_M$ for $\rho(t)$.
\\

From this upper bound, we obtain the lower bound on $S(t)$ because
$$
\beta \frac{d}{dt} S(t)  = Q\big(S(t),\rho(t) \big) \geq Q\big(S(t),\rho_M \big)
$$
and it is enough to use again \eqref{e.Qbounds} and the condition on the initial data \eqref{gen:id}.

\subsection{BV estimates}

Our next goal is to find a quantity which converges as $t \to \infty$. This step is crucial and we introduce a new idea which allows us to conclude.
\\

\noindent {\em 1st step. Equations on  $J := \dot{S}$ and $P := \dot{\rho}$.}
  With these definitions, from equations \re{e.n} and \re{e.I}, we can write 
\beq
P = \int nR\,dx, \qquad   \beta J  =  Q.  
\label{e.J}
\eeq
With the definitions 
\beq
0< \ul K_1 \rho(t) \leq  \alpha(t) := \int  n R_S\,dx \leq \ol K_1 \rho_M,\qquad\gamma(t) :=  \int  nR^2\,dx,
\label{lb:al}
\eeq
differentiating both equations on $n$ and $S$, we obtain 
\begin{align}
\dot{P} &= J \int  nR_S\,dx + \int  n_t R\,dx = \alpha(t) J  +  \gamma(t)  \label{e.Qdot},
\end{align}
\begin{align}
&\beta \dot{J}  =  Q_SJ+Q_\rho P.  \label{e.Jdot}
\end{align}
 
\noindent {\em 2nd step. Bound on a linear combination of $P$ and $J$.}
 Now we consider a linear combination of $P$ and $J$, where $\mu(t)$ is a function to be determined later. We write
\begin{align}
\frac{d}{dt} (P +\beta \mu J) &= \alpha J  + \beta \dot{\mu} J + \mu (  Q_S J+ Q_\rho P) +\gamma\quad \nonumber  \\
&= \mu  Q_\rho  (P + \beta \mu J) + (\beta \dot{\mu} - \beta Q_\rho \mu^2 + \mu  Q_S + \alpha) J  +\gamma. \label{e.qj}
\end{align}
We choose a function $\mu(t)$ such that the second parenthesis in the above equation is zero. In other words,  $\mu(t)$ solves the differential equation 
\beq 
\beta \dot{\mu} =- \beta |Q_\rho| \mu^2 + \mu  |Q_S| - \alpha  . 
\label{e.ode}
\eeq
Because the solution might blow-up to $-\infty$ in finite time, we first check that we can find a solution  $\mu(t)>0$ of \eqref{e.ode} for all times. To do so, we notice that the zeroes of 
$- \beta |Q_\rho| \mu^2 + \mu  |Q_S| - \alpha$ are 
$$
\mu_{\pm}(t) := \frac{1}{2 \beta |Q_\rho|} \big(  |Q_S| \pm \sqrt{Q_S^2 - 4 \alpha \beta  |Q_\rho|} \big).
$$
With assumptions \re{e.Qbounds} and \re{e.beta_rho},  both zeros are real positive. 

We are going to find two constants $0 <  \mu_m <  \mu_M$ such that, choosing initially $ \mu_m < \mu(0) < \mu_M$, 
then we have for all times
\beq
0 < \mu_m \leq \mu(t) \leq \mu_M := \f 1 \beta 
\max_{\begin{array}{l} \rho_m \leq \rho \leq \rho_M, \\ S_m \leq S \leq S_0 \end{array}}  \f{|Q_S|}{ |Q_\rho|},
\label{lb:mu}
\eeq
and $\mu_m$ defined by the condition 
\beq \label{e.muJ_cond}
\max_t \mu_- (t) \leq   \mu_m  := \min_t \mu_+(t).
\eeq

We first show how to enforce the inequality \eqref{e.muJ_cond}. We use that, for $0 \leq x \leq 1$, the concavity inequality holds: $\sqrt{1-x} \geq 1-x$, and compute
$$
\mu_-(t)  \leq 2 \f{\alpha (t)}{|Q_S|}  \leq 2\f{ \ol K_1 \rho_M}{|Q_S|},
$$
$$
\mu_+(t)  \geq \f{|Q_S|}{ \beta |Q_\rho|}  \left(1-2 \f{\alpha \beta  |Q_\rho|}{ |Q_S|^2} \right) \geq  \f{|Q_S|}{ \beta |Q_\rho|} - 2\f{ \ol K_1 \rho_M }{ |Q_S| } .
$$
The condition  \re{e.beta_rho} is enough to obtain the inequality \eqref{e.muJ_cond}.
\\

The lower bound in \eqref{lb:mu},  is because the condition \eqref{e.muJ_cond} imposes $\mu_m \in ( \mu_- (t) ,  \mu_+ (t))$ and thus $ \beta |Q_\rho| \mu_m^2 + \mu_m  |Q_S| - \alpha \geq 0$ for all $t\geq 0$. 
\\

For the upper bound in \eqref{lb:mu}, we use the inequality $\sqrt{1-x} \le 1 - \f x 2$ and we obtain
$$
\mu_+ < \f{|Q_S|}{2 \be |Q_\rho|}\left(2-2\al \be \f{|Q_\rho|}{|Q_S|^2}\right) \le \f{|Q_S|}{ \be |Q_\rho|} \leq \mu_M.
$$

With this choice of $\mu(t)$ and coming back to equation \re{e.qj}, we arrive to
$$
\frac{d}{dt} (P + \beta \mu J) \geq - \mu  \; |Q_\rho| \; (P + \beta \mu J),
$$
and we conclude that  
\begin{equation} \label{e.QmuJ}
\big(P(t) + \beta \mu(t) J(t) \big)_- \leq \big( P(0) + \beta \mu(0) J(0) \big)_-  e^{- K_Q \mu_m t},  \qquad \forall t\geq 0 .
\end{equation}
\\ 
\noindent {\em 3rd step. $L^1$-bound on $P$.}
From the above inequality, we wish to prove that $P(t)$ is integrable on the half-line. Adding $\al\f{P}{\be\mu}$ to \re{e.Qdot}, we find the ODE
$$
\frac{d}{dt} P +\al\f{P}{\be\mu}= \al\left(J+\f{P}{\be\mu}\right)+\gamma \geq -\al \left(J+\f{P}{\be\mu}\right)_- .
$$
Taking negative parts, we obtain the inequality
$$
\frac{d}{dt} P_- +\al\f{P_-}{\be\mu} \leq  \al\left(J+\f{P}{\be\mu}\right)_-,
$$
and, because $P$ is bounded, for some constant $C$
$$
\int_0^\infty \al(t) {P_- (t)}dt \leq C.
$$

With the lower bound on $\al$ in \eqref{lb:al}, we conclude that
$$
\ul K_1 \int  \rho \left(\ddt \rho\right)_-\,dx = \f{ \ul K_1}{2} \int \left(\ddt \rho^2\right)_-\,dx  \le \f C 2.
$$
and because $\rho(t)$ is bounded, we finally find that $ \ddt\rho^2$   is bounded in  $ L^1(0,\infty)$, therefore $\rho^2$ has a limit for $t \to \infty$ and $\rho$ has a limit $\bar \rho$
\\

\noindent {\em 4th step. Conclusion.}
Since $\rho(t)$ has a long term limit $\bar \rho$, the stability assumption for $Q$, more precisely  $Q_S<0$ in \eqref{e.Qbounds}, shows $S(t)$ also has 
a limit $\bar S$ and $Q(\bar \rho, \bar S)=0$.

As usual, \cite{OD.PJ.SM.BP:05, GB.BP:07, LD.PJ.SM.GR:08}, we can conclude that $R(x, \bar S) \leq 0$ for all $x$. Otherwise $n(t,x)$ would have an exponential growth for  $x$ in an open  set,  which would imply exponential growth for large times, a contradiction with the upper bound on $\rho(t)$.

This gives the statements of the Theorem~\ref{th:betasmall}.

\subsection{Numerical considerations}

For $\beta $ large, we could expect that the system could become unstable and that solutions  can be periodic. 
This is the case for inhibitory  integrate-and-fire models, these are pdes describing  neural networks, with strong relaxation properties to a steady state. It is well-known, see \cite{NB.VH.99} for instance, that delays can generate a spontaneous activity i.e. periodic solutions.

However, we did not observe such a behaviour in numerical simulations we conducted. This is confirmed by the stability analysis of a simplified equation. 

The numerics have been performed in Matlab with parameters as follows. As initial data we use
$S(t=0)=5$
and $n(t=0)=C_{\text{mass}} e^{-200(x-0.5)^2}$ where $C_{\text{mass}}$ is chosen such that the initial mass in the computational domain is equal to $5$. We set $ R:=20(-0.6 + 0.2S-(x-.5)^2)$ and $Q(\rho,S) := 8.5 -(0.5+\rho)S$.
The equation is solved  by an implicit-explicit finite-difference method on grid consisting of $1000$ points  (time step $dt= 4\cdot 10^{-4}$).
The plot shows oscillations of  $\rho(t)$ and $I(t)$. Moreover, numerically it seems that $\int_0^\infty |\dot \rho| \,dt$ is bounded.
\begin{figure}[h!]
\centerline{\includegraphics[width=.65\textwidth]{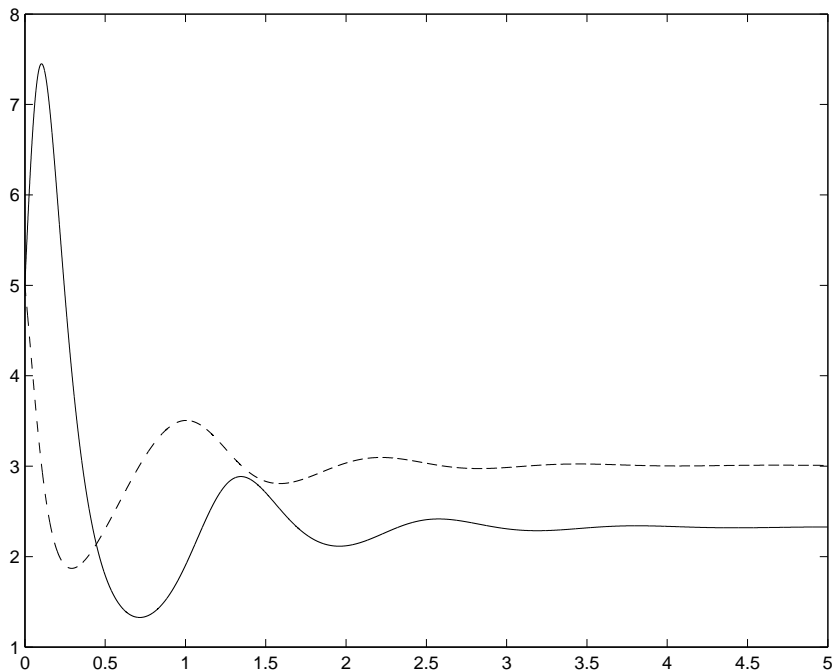}}
\caption{Dynamics of $\rho$ (\line(1,0){30}) and $I$  ($- - - -$).}
\label{F4}
\end{figure}

\begin{remark}
We can rewrite \re{e.Qdot} and \re{e.J} as 
\beq 
\dot{v} = A v + b
\eeq
where 
\beq
v= \left(
  \begin{array}{c}
  Q  \\
  J  
  \end{array} \right), \qquad A= \left(
  \begin{array}{cc}
   0 & -\al  \\
   1 & -1  
  \end{array} \right),
  \qquad
b =\left(
  \begin{array}{c}
   \int nR^2\,dx  \\
   0  
  \end{array} \right).
\eeq
For $\be$ small $A$ has real eigenvalues whereas for $\be$ large, it has complex eigenvalues.
Therefore our method cannot work for $\be$ large.
One direction to extend the result would be to work directly on the system \re{e.Qdot}--\re{e.J}.
\end{remark}

\section{Perspectives and open questions}

We have proved long term convergence to an ESD for a general model of a chemostat where the nurient delivery does not react immediately to the population dynamics. Our proof extends the proof based on Total Variation bounds developed in \cite{BP:07,GB.BP:07, GB.BP:08} and uses a fast (but not infinite) nutrient production measured by the small parameter $\beta$. 
\\

Surprisingly, the proof does not seem to give directly uniform TV bounds for  $\beta \approx 0$. It  does not seem to be possible with this approach to prove uniform bounds for the full range $\beta \in [0, \beta_0]$ for some small $\beta_0$, which could be a first step to prove uniform convergence of $S(t)$ for $t\in [0, \infty]$ as $\beta \to 0$.
\\

There are several related problems which, usually, can be approached with the same method. One of them is the rare mutations/long term behaviour described by the following extension of \eqref{e.n}--\eqref{e.rho}
$$ \left\{ \begin{array}{l}
 \e \partial_t n_\e(t,x) - \e^2 \Delta n_\e = n_\e \,R(x,S_\e(t)), \qquad x \in \R^d, \; t \geq 0,
\\ [5pt]
  \e \beta \frac{d}{dt} S_\e(t)  = Q\big(S_\e(t),\rho_\e(t) \big) ,  
\\ [5pt]
\rho_\e(t) :=  \int n_\e(t,x)\,dx,
\end{array} \right. 
$$
which can be treated using the constrained Hamilton-Jacobi approach \cite{OD.PJ.SM.BP:05,BP:07,GB.BP:07, GB.BP:08, AL.SM.BP:10}, provided some strong compactness is proved as e.g. TV bounds which are uniform in $\e$. 
\\

From the modelling side,  the TV bounds giving long term behaviour is not known in several examples of chemostat systems. An example is  the quasi-stationary case with general uptake rate and rendering factor,  
$$ \left\{ \begin{array}{l}
\p_t n= nR(x,S),
\\[5pt]
S+\int \eta(x,S)n(x)\,dx = S_0.
\end{array} \right. 
$$

\noindent {\bf Acknowledgment.} This research was supported by the french "ANR blanche" project Kibord:  ANR-13-BS01-0004. 

%
%
%


\bibliography{bibli_LP.bib}

\begin{thebibliography}{10}

\bibitem{ackleh_F_T}
Azmy~S. Ackleh, Ben~G. Fitzpatrick, and Horst~R. Thieme.
\newblock Rate distributions and survival of the fittest: a formulation on the
  space of measures.
\newblock {\em Discrete Contin. Dyn. Syst. Ser. B}, 5(4):917--928 (electronic),
  2005.

\bibitem{GB.BP:07}
G.~Barles and B.~Perthame.
\newblock Concentrations and constrained {H}amilton-{J}acobi equations arising
  in adaptive dynamics.
\newblock {\em Contemp. Math.}, 439:57--68, 2007.

\bibitem{NB.VH.99}
N.~Brunel and V.~Hakim.
\newblock Fast global oscillations in networks of integrate-and-fire neurons
  with low firing rates.
\newblock {\em Neural Comput.}, 11(7):1621--1671, October 1999.

\bibitem{NC.PJ.GR:11}
N.~Champagnat, P.-E. Jabin, and G.~Raoul.
\newblock Convergence to equilibrium in competitive lotka-volterra equations
  and chemostat systems.
\newblock {\em C. R. Acad. Sci. Paris S\'er. I Math.}, 348(23-24):1267--1272,
  2010.

\bibitem{RC.JH:05}
R.~Cressman and J.~Hofbauer.
\newblock Measure dynamics on a one-dimensional continuous trait space:
  theoretical foundations for adaptive dynamics.
\newblock {\em Th. Pop. Biol.}, 67(1):47--59, 2005.

\bibitem{LD.PJ.SM.GR:08}
L.~Desvillettes, P.-E. Jabin, S.~Mischler, and G.~Raoul.
\newblock On mutation-selection dynamics for continuous structured populations.
\newblock {\em Commun. Math. Sci.}, 6(3):729--747, 2008.

\bibitem{OD:04}
O.~Diekmann.
\newblock A beginner's guide to adaptive dynamics.
\newblock In {\em Mathematical modelling of population dynamics}, volume~63 of
  {\em Banach Center Publ.}, pages 47--86. Polish Acad. Sci., Warsaw, 2004.

\bibitem{OD.PJ.SM.BP:05}
O.~Diekmann, P.-E. Jabin, S.~Mischler, and B.~Perthame.
\newblock The dynamics of adaptation: an illuminating example and a
  {H}amilton-{J}acobi approach.
\newblock {\em Th. Pop. Biol.}, 67(4):257--271, 2005.

\bibitem{SG.EK.GM.JM:98}
S.~A.~H. Geritz, E.~Kisdi, G.~M\'eszena, and J.~A.~J. Metz.
\newblock Evolutionarily singular strategies and the adaptive growth and
  branching of the evolutionary tree.
\newblock {\em Evol. Ecol}, 12:35--57, 1998.

\bibitem{PJ.GR:09}
P.-E. Jabin and G.~Raoul.
\newblock On selection dynamics for competitive interactions.
\newblock {\em J. Math. Biol.}, 63(3):493--517, 2011.

\bibitem{AL.SM.BP:10}
A.~Lorz, S.~Mirrahimi, and B.~Perthame.
\newblock Dirac mass dynamics in multidimensional nonlocal parabolic equations.
\newblock {\em Comm. Partial Differential Equations}, 36(6):1071--1098, 2011.

\bibitem{SM.BP.JW:12}
S.~Mirrahimi, B.~Perthame, and J.~Y. Wakano.
\newblock Direct competition results from strong competition for limited
  resource.
\newblock {\em Preprint}.

\bibitem{BP:07}
B.~Perthame.
\newblock {\em Transport equations in biology}.
\newblock Frontiers in Mathematics. Birkh\"auser Verlag, Basel, 2007.

\bibitem{GB.BP:08}
B.~Perthame and G.~Barles.
\newblock Dirac concentrations in {L}otka-{V}olterra parabolic {PDE}s.
\newblock {\em Indiana Univ. Math. J.}, 57(7):3275--3301, 2008.

\bibitem{chemostatbook}
H.~L. Smith and P.~Waltman.
\newblock {\em The theory of the chemostat: dynamics of microbial competition}.
\newblock Cambridge Univ. Press, 1994.

\bibitem{RS.JE:02}
R.~W. Sterner and J.~J. Elser.
\newblock {\em Ecological Stoichiometry: The Biology of Elements from Molecules
  to the Biosphere}.
\newblock Princeton University Press, 2002.

\end{thebibliography}
\bibliographystyle{plain}
 

\end{document}